\documentclass[12pt, a4paper]{amsart}

\usepackage{amssymb}
\usepackage{mathrsfs}
\usepackage{stmaryrd}
\usepackage[dvips]{graphics}

\theoremstyle{plain}
\newtheorem{theorem}{Theorem}[section]
\newtheorem{lemma}[theorem]{Lemma}
\newtheorem{corollary}[theorem]{Corollary}
\newtheorem{proposition}[theorem]{Proposition}

\theoremstyle{definition}

\theoremstyle{remark}
\newtheorem*{claim}{Claim}
\newtheorem*{acknowledgements}{Acknowledgements}

\numberwithin{figure}{section}

\DeclareMathOperator{\trd}{tr\,deg}
\DeclareMathOperator{\aut}{Aut}
\DeclareMathOperator{\acl}{acl}
\DeclareMathOperator{\aclk}{acl_{\it K}}
\DeclareMathOperator{\aclkp}{acl_{\it K'}}
\DeclareMathOperator{\Endk}{End_{\it K}}

\DeclareMathOperator{\Homk}{Hom_{\it K}}

\newcommand{\C}{\mathfrak{C}}
\newcommand{\G}{\mathrm{G}}					
\newcommand{\GG}{\mathbb{G}}					
\newcommand{\PP}{\mathbb{P}}					
\newcommand{\CC}{\mathbb{C}}					
\newcommand{\RR}{\mathbb{R}}					
\newcommand{\Q}{\mathbb{Q}}					
\newcommand{\LL}{\mathscr{Q}}					
\newcommand{\J}{\mathscr{J}}					
\newcommand{\N}{\mathbb{N}}

\title{Combinatorial geometries of field extensions}
\author{Jakub Gismatullin}
\address{Mathematical Institute, University of Wroclaw, pl. Grunwaldzki 2/4, 50-384 Wroclaw, Poland}
\email{gismat@math.uni.wroc.pl}
\keywords{Combinatorial geometry, full algebraic matroid, projective plane, transcendental field extension.}
\subjclass[2000]{Primary 03C98, 51D20; Secondary 12F20, 05B35.}

\begin{document}

\begin{abstract}
We classify the algebraic combinatorial geometries of arbitrary field extensions of transcendence degree greater than $4$ and describe their groups of automorphisms. Our results and proofs extend similar results and proofs by Evans and Hrushovski in the case of algebraically closed fields. The classification of projective planes in algebraic combinatorial geometries in arbitrary fields of characteristic zero will also be given.
\end{abstract}

\maketitle

\section*{Introduction}

Let $K \subset L$ be an arbitrary field extension. We investigate the algebraic combinatorial geometry $\GG(L/K)$ and pregeometry $\G(L/K)$ in $L$ obtained from algebraic dependence relation over $K$. Such a geometry is sometimes called a full algebraic matroid.

In \cite{EH1} the authors classify projective planes in $\GG(L/K)$ for algebraically closed $K$ and $L$. Using their results, we give such a classification for arbitrary fields  $K$ and $L$ of characteristic zero. We prove a theorem about formulas with one quantifier of the first-order theory of $\GG(L/K)$. Assume that the transcendence degree of $L$ over $K$ is at least $5$. When considering $\GG(L/K)$ we may also assume that $L$ is a perfect field and $K$ is relatively algebraically closed in $L$. One of the main results of \cite{EH2} is the reconstruction of the field $L$ from $\GG(L/K)$ when $L$ is algebraically closed. We generalize this reconstruction to arbitrary field extension $K \subset L$ (of transcendence degree $\geq5$), and thus we obtain full classification of combinatorial geometries of fields:  $\GG(L_1/K_1)$ and $\GG(L_2/K_2)$ are isomorphic  if and only if the field extensions \[K_1\  \subset \ L_1\quad \text{and} \quad K_2\ \subset\ L_2\] 
are isomorphic (here we assume that $L_1$ and $L_2$ are perfect and $K_1$, $K_2$ are relatively algebraically closed). We also give a description of $\aut(\GG(L/K))$.

We work within a large algebraically closed field $\C$. By $\widehat{F}$ and $\widehat{F}^r$ we denote algebraic and purely inseparable closure of a field $F$ in $\C$. Throughout this paper we assume that $K \subset L$ is an arbitrary field extension and the transcendence degree of $L$ over $K$ is at least $3$. We take the basic definitions of algebraic combinatorial geometry and pregeometry from \cite{EH1,EH2}. For $X\subseteq L$, let  $\aclk(X)$ be $\widehat{K(X)}$. We denote by $\G(L/K)$ \emph{the pregeometry} $(L,\aclk)$. \emph{The geometry} $\GG(L/K)$ is obtained from $L\setminus \widehat{K}$ by factoring out the equivalence relation: \[x \sim y\quad\Longleftrightarrow\quad\widehat{K(x)}=\widehat{K(y)}.\] We can also transfer the closure operation $\aclk$ from $\G(L/K)$ to $\GG(L/K)$: \[\aclk(Y/{\sim})=\aclk(Y)/{\sim}.\] Therefore we can regard the points of $\GG(L/K)$ as sets $\aclk(x)$, where $x\in L\setminus \widehat{K}$, and $\aclk$ as the usual algebraic closure. When considering $\GG(L/K)$ we assume that $L$ is a perfect field and $K$ is relatively algebraically closed in $L$ (because $\GG(L/K)=\GG(\widehat{L}^r/K)=\GG(\widehat{L}^r/\widehat{L}^r\cap \widehat{K})$). Subsets of $\GG(L/K)$ of the form $\aclk(X)$, $X\subseteq\GG(L/K)$, are called \emph{closed}. \emph{The rank} of a subset of $\GG(L/K)$ or $\G(L/K)$ is its transcendence degree over $K$. We also have notions of \emph{independent set} (for each $x\in X$,\! $x\not\in \aclk(X\setminus\{x\})$) and a \emph{basis} of a closed subset as a maximal independent set (transcendence basis). Note that the closure operation $\aclk$ satisfies \emph{the exchange condition}: \[x\in \aclk(A\cup\{y\})\setminus\aclk(A)\quad\Longrightarrow\quad y\in\aclk(A\cup\{x\}).\] A closed subset of rank 1 (respectively 2, 3) is a \emph{point} (respectively \emph{line} and \emph{plane}). If $X$ is a closed subset of $\GG(L/K)$, and a tuple $\overline{x}\subset L$ satisfies $X=\aclk(\overline{x})$, then we say that $\overline{x}$ is \emph{generic} in $X$.

Let $F$ be a skew field (division ring). We will denote by $\PP(F)$ the projective plane over $F$. It is simply the set $F^3\setminus\{0\}$ factored out by the relation: \[(x_1,x_2,x_3) \simeq (y_1,y_2,y_3)\quad \Longleftrightarrow\quad \bigl(\exists\  0\neq\lambda\in F\bigr)\; (x_1,x_2,x_3)=\lambda(y_1,y_2,y_3).\]

The paper is organized as follows. The first section is devoted to giving some preliminary definitions and results from \cite{EH1}. In the second section we classify the projective planes arising in $\GG(L/K)$. Section 3 contains a theorem about first-order theory of $\G(L/K)$ and formulas with one quantifier. In Section 4 we transfer theorems from \cite{EH2} to geometries of arbitrary field extensions and prove a general classification theorem for them.

The reader is referred to \cite{Ma} for the model-theoretic background and notation, and to \cite{We} for general background on pregeometries and matroids.

\section{Preliminaries}
\label{s:prel}

For definitions and proofs in this section we refer the reader to \cite{EH1}. Throughout this section we assume that $K$ and $L$ are algebraically closed. Let $X$ be a subset of $\GG(L/K)$ and let $\aclk^X$ be the relative closure operation: $\aclk^X(Y)=\aclk(Y)\cap X$ for $Y\subseteq X$. We say that $X$ is a \emph{projective plane of} $\GG(L/K)$ if the geometry $(X,\aclk^X)$ is itself a projective plane, meaning that:

1) the geometry $(X,\aclk^X)$ has rank 3;

2) there are three noncollinear points in $X$;

3) any line has at least three different points;

4) any two lines intersect.

If a projective plane $X$ is isomorphic to $\PP(F)$, for some skew field $F$, then we say that $X$ is \emph{coordinatised} by $F$. It is well known (\cite[Chapter 7]{Ha}) that if the Desargues theorem is true in $X$, then $X$ is coordinatised by a unique skew field. The converse is also true. If $X_1$ and $X_2$ are Desarguesian projective planes coordinatised by $F_1$ and $F_2$ respectively, and $X_1\subseteq X_2$, then $F_1$ is a subskewfield of $F_2$. It is proved in \cite{Li} that any projective plane in $\GG(L/K)$ is Desarguesian. The aim of the next section is to find all skew fields coordinatising some projective planes in $\GG(L/K)$ for arbitrary $K\subset L$ of characteristic zero. The paper \cite{EH1} describes all such skew fields in the case when $L$ and $K$ are algebraically closed.

Let $(G,\ast)$ be a one-dimensional irreducible $K$-definable algebraic group in $L$. Then $G$ is isomorphic over $K$ (\cite[Section 3.1]{EH1}), as an algebraic group, to one of the following commutative groups: $(L,+)$, $(L^*,\cdot)$ or an elliptic curve. Since $G$ is commutative, the group $\Endk(G)=\Homk(G,G)$ of definable over $K$ morphisms of $G$ (as an algebraic group) may be given a ring structure $(\Endk(G),+,\circ)$ and is embeddable into a skew field of quotients $\Endk(G)_0$. If $char(L)>0$, then $\Endk(L,+)$ is the ring of $p$-polynomials over $K$ and we donote by $\mathcal{O}_{\widehat{K}}$ the skew field $\Endk(L,+)_0$. Let $\overline{x},\overline{y},\overline{z}\in G$ be independent generics over $K$. We may consider $G$ as an $\Endk(G)$ module and define \[\PP((G,\ast)\colon\overline{x},\overline{y}, \overline{z})=\{\aclk(a(\overline{x})\ast b(\overline{y})\ast c(\overline{z}))\colon (a,b,c)\in \Endk(G)^3\setminus\{\bf{0}\}\}.\] This is a projective plane in $\GG(L/K)$, coordinatised by $\Endk(G)_0$ i.e. elements of $\PP((G,\ast)\colon\overline{x},\overline{y}, \overline{z})$ are dependent with respect to $\Endk(G)$ exactly if they are $\aclk$-dependent.

\begin{lemma}\label{lem:impor} Let $x_1$, $x_2$, $x_3$, $x'_1$, $x'_2$, $x'_3\in L$.

(i) If each triple $\{x_1,x_2,x_3\}$ and $\{x'_1,x'_2,x'_3\}$ is algebraically independent over $K$, and
\begin{align*}
\aclk(x_1+x_2)&=\aclk(x'_1+x'_2),\quad \aclk(x_i)=\aclk(x'_i),\; \text{for } i=1,2,3,\\
\aclk(x_1+x_3)&=\aclk(x'_1+x'_3),
\end{align*}
then there exist $0\neq c,c'\in \Endk(L,+)$ and $d_1,d_2,d_3\in K$ such that \[c'(x'_1)=c(x_1)+d_1,\ c'(x'_2)=c(x_2) +d_2,\ c'(x'_3)=c(x_3) +d_3.\]

(ii) If each pair $\{x_1,x_2\}$ and $\{x'_1,x'_2\}$ is algebraically independent over $K$, and
\[ \aclk(x_1)=\aclk(x'_1),\quad \aclk(x_2)=\aclk(x'_2), \quad \aclk(x_1\cdot x'_1)=\aclk(x_2\cdot x'_2),\]
then there exist $0\neq n,m\in \mathbb{Z}$ and $0\neq a,b\in K$ such that $x_1^n=ax_2^m, \quad y_1^n=by_2^m.$
\end{lemma}
\begin{proof} First statement follows from \cite[Theorem 2.2.2]{EH1} and the second from \cite[Theorem 1.1]{EH2}.
\end{proof}

\section{Projective planes in $\GG(L/K)$}
\label{s:proj}

Throughout this section we assume that $K\subset L$ is an arbitrary field extension and $\trd_K(L)\geq3$. The geometry $\GG(L/K)$ naturally embeds into $\GG(\widehat{L}/\widehat{K})$. Therefore we can use theorems about $\GG(\widehat{L}/\widehat{K})$ to investigate $\GG(L/K)$. 

From the proof of \cite[Theorem 3.3.1]{EH1} we obtain some maximal projective planes in $\GG(\widehat{L}/\widehat{K})$ in the following way. Suppose $x,y,z\in \widehat{L}$ are algebraically independent over $K$. Then the projective plane \[\PP((\widehat{L},+)\colon x,y,z)\] is the largest projective plane in $\GG(\widehat{L}/\widehat{K})$ containing the tuple $\bigl(\aclk(x)$,\! $\aclk(y)$,\! $\aclk(z)$,\! $\aclk(x+y)$,\! $\aclk(x+z)\bigr)$, i.e. if a projective plane $\PP\subset\GG(\widehat{L}/\widehat{K})$ contains points $\bigl(\aclk(x)$,\! $\aclk(y)$,\! $\aclk(z)$,\! $\aclk(x+y)$,\! $\aclk(x+z)\bigr)$, then $\PP\subseteq \PP((\widehat{L},+)\colon x,y,z)$.

The next theorem generalizes above remark to the geometry $\GG(L/K)$ in characteristic zero. The case of positive characteristic requires detailed knowledge of the structure of $\mathcal{O}_{\widehat{K}}$.

\begin{theorem}($char(L)=0$) \label{th:4} Suppose that $x,y,z\in\widehat{L}$ are independent over $\widehat{K}$ and the tuple $\bigl(\aclk(x)$,\! $\aclk(y)$,\! $\aclk(z)$,\! $\aclk(x+y)$,\! $\aclk(x+z)\bigr)$ is in $\GG(L/K)$ ($x$, $y$ and $z$ do not need to be in $L$). Then
\[\PP((\widehat{L},+)\colon x,y,z)\cap \GG(L/K)=\{\aclk(ax+by+cz)\colon (a,b,c)\in(\widehat{K}\cap L)^{3}\setminus \{\bf 0\}\},\] is a projective plane in $\GG(L/K)$, coordinatised by $\widehat{K}\cap L$. Moreover the above plane is the largest projective plane in $\GG(L/K)$ containing the tuple $\bigl(\aclk(x)$,\! $\aclk(y)$,\! $\aclk(z)$,\! $\aclk(x+y)$,\! $\aclk(x+z)\bigr)$.

\end{theorem}
\begin{proof} Let $f\in\aut(\widehat{L}/L)$ be arbitrary. By assumption we have
\begin{align*}
\aclk(x)&=f[\aclk(x)]=\aclk(f(x)),\\
\aclk(y)&=\aclk(f(y)),\quad \aclk(x+y)=\aclk(f(x)+f(y)),\\
\aclk(z)&=\aclk(f(z)),\quad \aclk(x+z)=\aclk(f(x)+f(z)).
\end{align*}
Therefore by Lemma \ref{lem:impor} we obtain $f(x)=c'\cdot x + d_1$,\! $f(y)=c'\cdot y + d_2$,\! $f(z)=c'\cdot z + d_3$, for some $d_1,d_2,d_3\in \widehat{K}$ and $0\neq c'\in \widehat{K}$.\\

$\subseteq$: Let $v\in \PP((\widehat{L},+)\colon x,y,z)\cap \GG(L/K)$. We have $v=\aclk(ax+by+cz)=\aclk(l)$, where $a,b,c\in \widehat{K}$ and $l\in L$. It follows $f[v]=v$, so
\begin{gather*}
\aclk(ax+by+cz)=\aclk (f(a)f(x)+f(b)f(y)+f(c)f(z))\\
=\aclk(c'\cdot(f(a)x+f(b)y+f(c)z)+d')=\aclk(f(a)x+f(b)y+f(c)z),
\end{gather*}
for $c'$,\! $d'=f(a)d_1+f(b)d_2+f(c)d_3\in \widehat{K}$ (because $f[\widehat{K}]=\widehat{K}$). By \cite[Example 2, Section 3.3]{EH1} there is a nonzero $\lambda\in\widehat{K}$ such that $(f(a),f(b),f(c))=\lambda(a,b,c)$. If e.g. $a\neq0$, then $f(\tfrac{b}{a})=\tfrac{b}{a}$ and $f(\tfrac{c}{a})=\tfrac{c}{a}$. But $f$ has been arbitrary, so $\tfrac{b}{a},\tfrac{c}{a}\in L$. Finally $v=\aclk(ax+by+cz)=\aclk(x+\tfrac{b}{a}y+\tfrac{c}{a}z)$, where $\tfrac{b}{a}, \tfrac{c}{a}\in\widehat{K}\cap L$.\\

$\supseteq$: Let $a,b,c\in \widehat{K}\cap L$ and consider $v=\aclk(ax+by+cz)$. It remains to prove that $v\in\GG(L/K)$. We have $f[v]=v$, because
\begin{gather*}
f[\aclk(ax+by+cz)]=\aclk(af(x)+bf(y)+cf(z))\\
=\aclk(c'\cdot(ax+by+cz)+d')=\aclk (ax+by+cz).
\end{gather*}
Let $w(x)$ be a minimal monic polynomial for $ax+by+cz$ over $L$. Then $v=\aclk(ax+by+cz)=\aclk(\text{roots of }w)= \aclk(\text{coefficients of }w)\in\GG(L/K)$.

The last part of the theorem follows from the first part and from remarks at the begining of this section.
\end{proof}

From the above we have that the geometries $\GG(\CC/\Q)$ and $\GG(\RR/\Q)$ are not isomorphic, because in $\GG(\RR/\Q)$ there is a maximal projective plane, coordinatised by $\widehat{\Q}\cap\RR$ and in $\GG(\CC/\Q)$ there is no such plane.

The next result generalizes \cite[Corollary 3.3.2]{EH1} and follows from Theorem \ref{th:4}.

\begin{corollary} ($char(L)=0$) If $\PP\subset\GG(L/K)$ is a projective plane, then $\PP$ is coordinatised by a subfield of one of the following fields: $\Q(\sqrt{-d})$, $d\in\N$ and $\widehat{K}\cap L$.
\end{corollary}

\section{The first-order theory of $\G(L/K)$}
\label{s:theory}

We can regard $\G(L/K)$ (and thus $\GG(L/K)$) as a model in the countable first-order language $\mathcal{L}=\{\acl_n\colon n\in\N\}$. Namely let \[\acl_n(a_0,\dotsc,a_n) \iff a_0 \in \aclk(a_1,\dotsc,a_n).\] We obtain a structure $(L, \mathcal{L})$. The following Theorem \ref{th:my} describes a small part of the first-order theory of $(L, \mathcal{L})$.

\begin{proposition} \label{prop:2} Let $F$ be an arbitrary field. If $F=F_1\cup\dotsb\cup F_n$, for some subfields $F_1,\dotsc,F_n$ of $F$, then $F=F_i$ for some $1\leq i \leq n$.
\end{proposition}
\begin{proof} This follows from a well known result of B. H. Neumann \cite{Ne}: if there is a covering of an abelian group by finitely many cosets of subgroups, then one of these subgroup has finite index. We leave the proof to the reader.
\end{proof}

\begin{theorem} \label{th:my} Let $K\subseteq L_1 \subseteq L_2$ be arbitrary field extensions. Assume that $\trd_K L_1 = \trd_K L_2 <\aleph_0$ or $\trd_K L_1, \trd_K L_2 \geq \aleph_0$. Then \[(L_1, \mathcal{L})\prec_1 (L_2, \mathcal{L}),\] i.e. for every $\mathcal{L}$-statement $\psi \in \mathcal{L}(L_1)$ with one quantifier and parameters from $L_1$ we have $(L_1, \mathcal{L})\models \psi \iff (L_2, \mathcal{L})\models \psi$.
\end{theorem}

It is easy to check that without the condition on transcendence degree, the theorem will not be true.

\begin{proof} We can assume that $\psi = \exists x \varphi(x,\overline{c})$, where $\overline{c} \subseteq L_1$ and $\varphi$ is a quantifier free formula. Using the exchange property for $\aclk$ we may assume that $\varphi(x,\overline{c})$ is of the form
\[ \bigvee_{k<l} \biggl(\Bigl(x \in \bigcap_{i<n_k} \acl_{K}(A_{k,i})\setminus\bigcup_{j<m_k} \acl_{K}(B_{k,j})\Bigr)\wedge \bigl(\text{(in)equality about }x,\, \overline{c}\bigl) \biggl),\]
where $A_{k,i}, B_{k,j} \subseteq \overline{c}\subseteq L_1$.

Let $\{p,q\}=\{1,2\}$ and assume that $L_p\models \exists x \varphi(x,\overline{c})$, then there exists $a\in L_p$ with $L_p\models \varphi(a,\overline{c})$. Without loss of generality we may assume that $a \notin \overline{c}$, so \[a \in \bigcap_{i<n} \acl_{K}(A_i)\setminus\bigcup_{j<m} \acl_{K}(B_j).\] If $n=0$, then by assumptions we have $a'\in L_q$ such that $L_q\models \varphi(a',\overline{c})$. Let $n\neq 0$. Note that by Proposition \ref{prop:2}, $L_p\models \exists x \varphi(x,\overline{c})$ is equivalent to: \[\bigl(\forall j<m\bigr)  \acl_{K}(B_j)\cap\bigcap_{i<n} \acl_{K}(A_i) \varsubsetneq \bigcap_{i<n} \acl_{K}(A_i),\] i.e. for $j<m$,\! $\widehat{K(B_j)}\cap\bigcap_{i<n} \widehat{K(A_i)}\cap L_p \varsubsetneq \bigcap_{i<n} \widehat{K(A_i)}\cap L_p$. The next lemma will be useful in the proof.

\begin{lemma} \label{lem:1} Suppose that $A$ and $B$ are finite subsets of $L_1$. Then there exists a finite subset $C\subseteq L_1$ satisfying \[\widehat{K(A)}\cap \widehat{K(B)} = \widehat{K(C)}.\]
\end{lemma}
\begin{proof} Let $C'$ be a transcendence basis of a $\widehat{K(A)}\cap \widehat{K(B)}$ over $K$. Write $C'=\{c_1,\dotsc,c_k\}\subset \widehat{L_1}$. Then $\widehat{K(A)}\cap \widehat{K(B)} = \widehat{K(C')}$. Take a minimal monic polynomial $w_i\in L_1[X]$ for $c_i$ over $L_1$ and let $C=\bigcup_{1\leq i \leq k}(\text{coefficients of }w_i)\subseteq L_1$. We shall show that $\widehat{K(C')} = \widehat{K(C)}$. By definition $C'\subset \widehat{K(C)}$, hence $\subseteq$. By symmetric polynomials we obtain \[C\subseteq K\biggl(\bigcup_{1\leq i \leq k}\text{roots of }w_i\biggr)\] and the last field is included in $\widehat{K(C')}$. We explain the last inclusion: if $w_i(a)=0$, then there exist  $f\in \aut(\widehat{L_1}/L_1), f(c_i)=a$, and thus $c_i\in \widehat{K(C')}$. Finally $a=f(c_i)\in f[\widehat{K(C')}]=f[\widehat{K(A)}\cap \widehat{K(B)}] \overset{A,B\subset L_1}{=}\widehat{K(A)}\cap \widehat{K(B)} = \widehat{K(C')}$.
\end{proof}

By Lemma \ref{lem:1}, to finish the proof it remains to show the following lemma.

\begin{lemma} For $A,B\subseteq L_1$ \[\widehat{K(A)}\cap L_1 \varsubsetneq \widehat{K(B)}\cap L_1 \iff \widehat{K(A)}\cap L_2 \varsubsetneq \widehat{K(B)}\cap L_2.\]
\end{lemma}
\begin{proof} Implication $\Rightarrow$ is obvious. $\Leftarrow$: Suppose, contrary to our claim, that $\widehat{K(A)}\cap L_1 = \widehat{K(B)}\cap L_1$. Take $a\in (\widehat{K(B)} \setminus \widehat{K(A)})\cap L_2$ and the minimal monic polynomial $w(X)=X^n+a_{n-1}X^{n-1}+\dotsb+a_1X+a_0\in K(B)[X]$ for $a$ over $K(B)$. Then for $i<n$ we have $a_i\in K(B)\subseteq L_1$, so by assumption $a_i\in \widehat{K(A)}\cap L_1$, hence $w\in (\widehat{K(A)}\cap L_1)[X]$ and thus $a$ is in algebraic closure of $\widehat{K(A)}\cap L_1$, which is included in $\widehat{K(A)}$. However by assumption $a\notin \widehat{K(A)}\cap L_2$ and $a\in L_2$, which is impossible.
\end{proof}
\end{proof}

\section{The reconstruction $L$ from $\GG(L/K)$ and corollaries}
\label{s:inter}

In this section we generalize some theorems of \cite{EH2} from the case of algebraically closed fields to the case of arbitrary field extensions. Throughout this section $K\subset L$ will be an arbitrary field extension, $char(L)=p$\ \  and $\trd_K(L)\geq5$. 

We begin with important definitions (see \cite{EH2} Definitions 2.1, 2.3 and 2.6). Let $(L\setminus K)^{(2)}$ denote the set of pairs $(x,y)\in (L\setminus \widehat{K})^2$ such that $x$ and $y$ are algebraically independent over $K$. We define the following subsets of $\GG(L/K)^4$:
\begin{align*}
&\LL=\{(\aclk(x),\aclk(y),\aclk(x+y),\aclk(x/y))\colon (x,y)\in(L\setminus K)^{(2)}\}\\
     &=\{(\aclk(x),\aclk(xz),\aclk(xz+x),\aclk(z))\colon (x,z)\in(L\setminus K)^{(2)}\},\\
&\LL'=\{(\aclk(x),\aclk(y),\aclk(x+y),\aclk(x\cdot y))\colon (x,y)\in(L\setminus K)^{(2)}\},\\
   &\J=\text{Im}(j),
\end{align*}
where $j\colon(L\setminus K)^{(2)}\rightarrow \GG(L/K)^5$ is the function \[j(x,a)=(\aclk(x),\aclk(x+a),\aclk(xa),\aclk(x+xa),\aclk(a)).\]

Let $\psi(A_1,A_2,B_1,B_2,C_1,C_2,D,E,F,G,H,I,P,Q,R,S,T,U,X,Y,Z)$ be an $\mathcal{L}$-formula (see Section \ref{s:theory}), standing for the assumptions from \cite[Lemma 3.1 (1), (2)]{EH2} and \cite[Corollary 3.4]{EH2}, i.e. let $\psi(A_1,\ldots,Z)$ means the conjunction of the following conditions (where $A=\aclk(A_1,A_2), B=\aclk(B_1,B_2), C=\aclk(C_1,C_2)$):

\begin{itemize}
\item[$\boxcircle$] $\trd(ABC)=\trd(AB)=\trd(BC)=\trd(AC)=4$,
\item[$\boxcircle$] $X \in\aclk(AY)\  \wedge\  Z\in\aclk(BY)\cap\aclk(CX)\  \wedge\  (X, Y, Z\not\in\aclk(ABC))$,
\item[$\boxcircle$] $(\forall_{A'\in A, B'\in B, C'\in C})\  X\not\in\aclk(A'Y)\ \wedge\  Z\not\in\aclk(B'Y)\cup\aclk(C'X)$,
\item[$\boxcircle$] $S\in A\ \wedge\  T\in B\ \wedge\ U\in C\ \wedge\ \trd(STU)=2$,
\item[$\boxcircle$] $(\exists_{S', T',U'})\ (S, T, U, S', T', U')$ is a partial quadrangle (see \cite[Section 2.1]{EH1}),
\item[$\boxcircle$] $\aclk(P,Y)\cap\aclk(S,X)=D\ \wedge\ \aclk(T,Y)\cap\aclk(Q,Z)=E$,
\item[$\boxcircle$] $\aclk(U,X)\cap\aclk(T,D)=F\ \wedge\ \aclk(P,T)\cap\aclk(Q,R)=G$, 
\item[$\boxcircle$] $\aclk(H,X)\cap\aclk(P,Y)=I\ \wedge\ \aclk(U,G)\cap A=H\ \wedge\ \aclk(F,Z)\cap C=R$.
\end{itemize}

\begin{theorem} The sets $\LL$, $\LL'$ and $\J$ are definable without parameters in $\GG(L/K)$.
\end{theorem}
\begin{proof} When $K$ and $L$ are algebraically closed then the proof of this theorem can be found in \cite[Section 3]{EH2}. We sketch it in this case.

Let $\psi_{\LL}(P,D,Y,I)=\bigl(\exists A_1,\dotsc,Z\bigr)\ \psi(A_1,\dotsc,Z)$ where the quantifier is free from $P$,\! $D$,\! $Y$,\! $I$. Then the formula $\psi_{\LL}$ defines $\LL$. Now we find a formula for $\LL'$. Lemma 2.2 in \cite{EH2} gives us configuration for multiplication: if $(x,y)\in(L\setminus K)^{(2)}$ and if the points $A'$,\!  $B'$,\! $C'$,\! $D'$,\! $E'\in\GG(L/K)$ are such that the configuration of points and lines in $\GG(L/K)$ holds as in Figure \ref{f:lem}, then \[E'=\aclk(x\cdot y),\] and there exist $a\in A'$ such that $A'=\aclk(a)$, $B'=\aclk(ax)$, $C'=\aclk(ay)$ and $D'=\aclk(axy)$.
\begin{figure}
\includegraphics{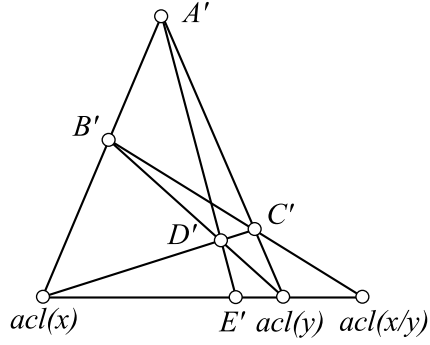}
\caption{Configuration for the multiplication}
\label{f:lem}
\end{figure}
Thus if we know $\aclk(x)$, $\aclk(y)$ and $\aclk(x/y)$, then in $\GG(L/K)$ we can construct $\aclk(xy)$. Let $\psi_{\LL'}(A,B,C,D)$ be
\[\bigl(\exists A',B',C',D',V\bigr)\ \psi_{\LL}(A,B,C,V) \wedge (\text{the configuration in Figure \ref{f:lem} holds}),\] where in Figure \ref{f:lem} we put $A$ instead of $\aclk(x)$, $D$ instead of $E'$, $B$ instead of $\aclk(y)$ and $V$ instead of $\aclk(x/y)$. Therefore $\psi_{\LL'}$ defines $\LL'$. To find a formula for $\J$ we recall \cite[Proposition 2.4]{EH2}: let $X$,\! $P$,\! $Q$,\! $R$,\! $A$ be in $\GG(L/K)$. Then \[(X,P,Q,R,A) \in \J \quad \Longleftrightarrow \quad (X,Q,R,A) \in \LL \wedge \bigl((X,A,P,Q),(X,A,P,R) \in\LL'\bigr). \] Hence, the formula \[\psi_{\J}(X,P,Q,R,A)= \psi_{\LL}(X,Q,R,A) \wedge \psi_{\LL'}(X,A,P,Q) \wedge \psi_{\LL'}(X,A,P,R)\] defines $\J$ in the algebraically closed case.

Now we turn to the general case, i.e. when $K$ and $L$ are arbitrary fields. It is sufficient to prove the next Claim, because we have for instance $(X,Q,R,A) \in \LL \Leftrightarrow \bigl(\exists P\bigr)(X,P,Q,R,A)\in\J$.

\begin{claim} The formula $\psi_{\J}$ defines $\J$ in $\GG(L/K)$.
\end{claim}

We will prove that the following conditions are equivalent:
\begin{enumerate}
\item $(X,P,Q,R,A)\in\J^{\GG(L/K)}$,
\item $(X,P,Q,R,A)\in\J^{\GG(\widehat{L}/\widehat{K})} \textrm{ and } (X,P,Q,R,A)\subset\GG(L/K)$,
\item $\GG(\widehat{L}/\widehat{K})\models\psi_{\J}(X,P,Q,R,A) \textrm{ and } (X,P,Q,R,A) \subset\GG(L/K)$,
\item $\GG(L/K)\models\psi_{\J}(X,P,Q,R,A)$.
\end{enumerate}

Implications (i)$\Rightarrow$(ii) and (ii)$\Leftrightarrow$(iii) are obvious. For (ii)$\Rightarrow$(i) take an arbitrary $f\in\aut(\widehat{L}/L)$ and write $(X,P,Q,R,A)=j(x,a)$ for some $(x,a)\in(\widehat{L}\setminus\widehat{K})^{(2)}$. Since $j(x,a)\subset\GG(L/K)$, we have $j(x,a)=f(j(x,a))=j(f(x),f(a))$, so by \cite[Lemma 2.5]{EH2} there exist $n\in\mathbb{Z}$ such that $f(x)=x^{p^n}$, $f(a)=a^{p^n}$. We show that $f(x)=x$ and $f(a)=a$. On the contrary, suppose that $p\neq0$ and $n\neq0$. Then $f(f(x))=f(x^{p^n})=f(x)^{p^n}=x^{p^{2n}}$ and in general $f^{m}(x)=x^{p^{m\cdot n}}$. However $x\in\widehat{L}$, so the set $\{f^{m}(x):m\in\N\}$ is finite. Hence, there is $k\in\N$ such that $x^k=1$, which implies $x\in\widehat{K}$, a contradiction.

(iv)$\Rightarrow$(iii): It is sufficient to show that for $(A,B,C,D)\subset\GG(L/K)$
\[\GG(L/K)\models\psi_{\LL}(A,B,C,D)\quad\Longrightarrow\quad \GG(\widehat{L}/\widehat{K}) \models\psi_{\LL}(A,B,C,D),\]
\[\GG(L/K)\models\psi_{\LL'}(A,B,C,D)\quad\Longrightarrow\quad \GG(\widehat{L}/\widehat{K}) \models\psi_{\LL'}(A,B,C,D).\]
It is immediately seen that we must only prove the following: for $X$,\! $Y$,\! $A_1$,\! $A_2\in\GG(L/K)$, if $\GG(L/K)\models \bigl(\forall A'\in \aclk(A_1,A_2)\bigr)\ X\not\in\aclk(A'Y)$ then $\GG(\widehat{L}/\widehat{K}) \models \bigl(\forall A'\in \aclk(A_1,A_2)\bigr)\ X\not\in\aclk(A'Y).$
Since $\GG(\widehat{L}/\widehat{K})=\GG(\widehat{L}/K)$, and the above formula has one quantifier, our statement follows from Theorem \ref{th:my}.

(i)$\wedge$(iii)$\Rightarrow$(iv): Take an element $(x',a')\in(L\setminus K)^{(2)}$ such that $(X,P,Q,R,A)=j(x',a')$. We must show the following (remember that $\aclk(a'+1)=\aclk(a')$, etc.)
\begin{itemize}
\item $\GG(L/K)\models\psi_{\LL}(\aclk(x'),\aclk(x'a'),\aclk(x'+x'a'),\aclk(x'a'/x'))$,
\item $\GG(L/K)\models\psi_{\LL'}(\aclk(x'),\aclk(a'),\aclk(x'+a'),\aclk(x'a'))$,
\item $\GG(L/K)\models\psi_{\LL'}(\aclk(x'),\aclk(a'+1),\aclk(x'+(a'+1)),\aclk(x'(a'+1)))$.
\end{itemize}
It is an easy consequence of \cite[Corollary 3.4]{EH2}. We give the proof only for the first case, the other cases are left to the reader. Let \[(P,D,Y,I) =(\aclk(x'),\aclk(x'a'),\aclk(x'+x'a'),\aclk(x'a'/x')=\aclk(a')).\]
We can find an algebraically independent (over $\widehat{K}$) set $\{a,b,c,d,x\}\in L$ such that \[(P,D,Y,I) =(\aclk(b),\aclk(ax),\aclk(ax+b),\aclk(ax/b)),\] and define the points $A_1, A_2,\ldots,R$ as in the conclusion of \cite[Corollary 3.4]{EH2} i.e. $A_1=\aclk(a), A_2=\aclk(b),\dotsc,R=\aclk(cb+d)$. Finally points $A,\dots,Z\in\GG(L/K)$ satisfy the assumption of \cite[Corollary 3.4]{EH2}, and thus they fulfil the formula $\psi$, so $P,D,Y$ and $I$ fulfil $\psi_{\LL}$ in $\GG(L/K)$.
\end{proof}

Now we prove the main classification theorem. We recall that $\widehat{F}^r=\bigcup_{n\in\N}F^{p^{-n}}$ is purely inseparable closure of $F$.

\begin{theorem} Suppose that $K\subset L$ and $K'\subset L'$ are field extensions and $\trd_K(L)$, $\trd_{K'}(L')\geq5$.
\begin{enumerate}
\item[(i)] The field $\widehat{L}^r$ is uniformly interpretable in $\GG(L/K)$, using a formula with one (arbitary) parameter from $L\setminus\widehat{K}$.

\item[(ii)] Every isomorphism $F\colon \GG(L/K)\xrightarrow{\cong} \GG(L'/K')$ is induced by some isomorphism $\widetilde{F}\colon \widehat{L}^r\xrightarrow{\cong} \widehat{L'}^r$ such that $\widetilde{F}\bigl[\widehat{L}^r\cap\widehat{K}\bigr]=\widehat{L'}^r\cap\widehat{K'}$, and for each $x\in \widehat{L}^r\setminus\widehat{K}$, $F(\aclk(x))=\aclk(\widetilde{F}(x))$. In particular $\GG(L/K) \cong \GG(L'/K')$ if and only if field extensions $\widehat{L}^r\cap \widehat{K}\subset \widehat{L}^r$ and $\widehat{L'}^r\cap\widehat{K'}\subset\widehat{L'}^r$ are isomorphic.

\item[(iii)] The natural mapping \[H\colon\aut(\widehat{L}^r/\{\widehat{L}^r\cap\widehat{K}\})\longrightarrow \aut(\GG(L/K)),\] is an epimorphism. If $char(L)=0$, then $H$ is an isomorphism of groups and if $char(L)\neq0$ then $\ker H\cong\mathbb{Z}$ is generated by the Frobenius automorphism.
\end{enumerate}
\end{theorem}

\begin{proof} Let $\cong$ be the following equivalence relation on $\J$ (\cite[Definition 2.9]{EH2})
\[j(x,a)\cong j(x',a')\quad \Longleftrightarrow \quad\bigl(\exists n\in\mathbb{Z}\bigr)\ a'=a^{p^n}.\]
Using \cite[Lemma 2.8]{EH2} we obtain that $\cong$ is a definable (without parameters) equivalence relation on $\J$. When $x$,\! $x'$ and $a$ are algebraically indepedent (over $K$), then \[j(x,a)\cong j(x',a')\ \Leftrightarrow\ (\exists P)\text{ configuration from Fig. 3 in \cite[Lem. 2.8]{EH2} holds}.\]
Implication $\Leftarrow$ follows from \cite[Lemma 2.8]{EH2}. For $\Rightarrow$ assume that $a'=a^{p^n}$ for some $0\leq n\in\mathbb{Z}$. We must find a suitable $P$ from $\GG(L/K)$. We have $\aclk(ax)=\aclk(a^{p^n}x^{p^n})=\aclk(a'x^{p^n})$, $\aclk((a+1)x)=\aclk((a^{p^n}+1)x^{p^n}) = \aclk((a'+1)x^{p^n})$ and $\aclk(x)=\aclk(x^{p^n})$. Hence $P=\aclk(x^{p^n}/x')\in\GG(L/K)$.

When $x$,\! $x'$ and $a$ are collinear, then we put (\cite[Definition 2.9]{EH2})
\begin{align*}
j(x,a)\cong j(x',a')\ \Leftrightarrow\ &\text{there exist }j(z,a'')\text{ with }z\not\in\aclk(a,x)\text{ and the configuration}\\
&\text{in Fig. 3 holds between }j(x,a),j(z,a'')\text{ and }j(z,a''),j(x',a').
\end{align*}

Take an arbitrary $a\in L\setminus\widehat{K}$ and let \[\J_1=[j(x,a)]_{\cong}=\{j(x',a):x'\in \widehat{L}^r\setminus\aclk(a)\}\] be one of the classes of $\cong$  (here we use the equality $j(x',a^{p^n})=j(x'^{p^{-n}},a)$ and properties of $\widehat{L}^r$). Let \[\mu\colon\J_1\rightarrow \widehat{L}^r\setminus\aclk(a),\quad\mu(j(x,a))=x.\] It follows from \cite[Lemma 2.5]{EH2}, that $\J_1$ be a bijective map and from \cite[Lemma 2.11]{EH2}, that there are generic definable over $a$ operations $\oplus$ and $\odot$ on $\J_1\times\J_1$ which satisfy: if $j(x,a), j(y,a)\in(L\setminus K)^{(2)}$ and $x,x',a$ are independent, then
\begin{align*}
& j(x,a)\oplus j(x',a)=j(x+x',a)\\
& j(x,a)\odot j(x',a)=j(x\cdot x',a).
\end{align*}
(the same definition works for non-algebraically closed case). Note that the map $\mu$ respects these operations (when defined). We now interpret the field $\widehat{L}^r$ in $\GG(L/K)$. Define a relation $\equiv$ on $\J_1^2$ by
\[(j(x_1,a),j(x_2,a))\equiv(j(y_1,a),j(y_2,a))\quad\Longleftrightarrow\quad \dfrac{x_1}{x_2}=\dfrac{y_1}{y_2}.\]
It is a definable over $a$ equivalence relation. We moreover define the product and the sum of two classes $[j(x_1,a),j(x_2,a)]_{\equiv}$ and $[j(x_1,a),j(x_2,a)]_{\equiv}$ as in \cite{EH2}, i.e.
\[[j(x',a),j(x,a)]_{\equiv}\ \cdot\ [j(y',a),j(y,a)]_{\equiv}\ =\ [j(x'',a),j(y'',a)]_{\equiv},\] for suitable $x''$ and $y''$ such that $\tfrac{x'}{x}\tfrac{y'}{y}=\tfrac{x''}{y''}$. We need a new class $0_{\equiv}$ to define the sum of classes in a standard fashion. Finally we extend $\mu$ to the isomorphism of fields: \[\mu\colon(\J_1^2/\equiv)\cup\{0_{\equiv}\}\xrightarrow{\cong}\widehat{L}^r,\quad \mu(j(x,a),j(x',a))=\dfrac{x}{x'},\ \mu(0_{\equiv})=0,\] which establishes (i).

(ii): We continue to use the notation of (i). Let $F\colon\GG(L/K)\xrightarrow{\cong}\GG(L'/K')$ be an isomorphisms of geometries (i.e. $F$ preserves the closure operator). Then
\[F\colon\J_1^{\GG(L/K)}=[j(x,a)]_{\cong}\  \overset{\cong}{\longrightarrow}\  \J_1^{\GG(L'/K')}=[j(y,b)]_{\cong},\]
\[F\colon\Bigl(\J_1^{\GG(L/K)}\Bigr)^2/\equiv\ \cup\ \{0_{\equiv}\}\ \overset{\cong}{\longrightarrow}\  \Bigl(\J_1^{\GG(L'/K')}\Bigr)^2/\equiv\ \cup\ \{0_{\equiv}\},\] for some $y,b\in L'$. Hence $F$ induces an isomorphism of fields $\widetilde{F}\colon \widehat{L}^r \xrightarrow{\cong} \widehat{L'}^r$.

Now, we show that $\widetilde{F}[\widehat{L}^r\cap\widehat{K}]=\widehat{L'}^r\cap\widehat{K'}$. Let $c\in\widehat{L}^r\cap\widehat{K}$ and $x\in\widehat{L}^r\setminus\aclk(a)$. Then we may write
\begin{align}
F(j(cx,a))&=j(y_1,b), \label{a:1}\\
F(j(x,a))&=j(y_2,b), \label{a:2}\\
\widetilde{F}(c)&=\mu(F([j(cx,a),j(x,a)]_{\equiv}))=\mu([j(y_1,b),j(y_2,b)]_{\equiv}) =\tfrac{y_1}{y_2},
\end{align}
for some $y_1,y_2\in \widehat{L'}^r$. However $c\in\widehat{K}$ yields that
\begin{align*}
j(cx,a)&=(\aclk(cx),\aclk(cx+a),\aclk(cxa),\aclk(cxa+cx),\aclk(a))\\
&=(\aclk(x),\aclk(cx+a),\aclk(xa),\aclk(xa+x),\aclk(a)).
\end{align*}
Thus from (\ref{a:1}) and (\ref{a:2}) above, we have
\begin{align*}
F(\aclk(x))&=\aclkp(y_1)=\aclkp(y_2),\\
F(\aclk(xa))&=\aclkp(y_1b)=\aclkp(y_2b),\\
F(\aclk(xa+x))&=\aclkp(y_1b+y_1)=\aclkp(y_2b+y_2).
\end{align*}
Hence from Lemma \ref{lem:impor} (ii) we obtain $n,m\in\mathbb{Z}$ and $d_1,d_2\in\widehat{K'}$ satisfying $y_1^n=d_1y_2^m$ and $(y_1b)^n=d_2(y_2b)^m$. This gives that $n=m$, and $y_1=c'y_2$ for some $c'\in\widehat{K'}$. Finally $\widetilde{F}(c)=\tfrac{y_1}{y_2}=c'\in\widehat{L'}^r\cap\widehat{K'}$.

Now we show the following
\[\bigl(\forall x\in \widehat{L}^r\setminus\widehat{K}\bigr)\quad F(\aclk(x))=\aclkp(\widetilde{F}(x)).\]
It follows from the preceding results that for $x_1,x_2\in\widehat{L}^r\setminus\aclk(a)$
\[F([j(x_1,a),j(x_2,a)]_{\equiv})\ =\ [j(y\widetilde{F}\Bigl(\dfrac{x_1}{x_2}\Bigr),b),j(y,b)]_{\equiv},\] for some $b\in\widehat{L'}^r$ and $y\in\widehat{L'}^r\setminus\aclkp(b)$. Let $t=\dfrac{y}{\widetilde{F}(x_2)}$. We obtain
\[\bigl(\forall x\in \widehat{L}^r\setminus\aclk(a)\bigr)\quad F(j(x,a))=j(\widetilde{F}(x)t,b).\]
Let $x_1,x_2\in\widehat{L}^r$ be algebraically independent over $\aclk(a)$. Then
\begin{align*}
j(\widetilde{F}(x_1x_2)t,b)&=F(j(x_1x_2,a))=F(j(x_1,a)\odot j(x_2,a))=F(j(x_1,a))\odot F(j(x_2,a))\\
&=j(\widetilde{F}(x_1)t,b)\odot j(\widetilde{F}(x_2)t,b)=j(\widetilde{F}(x_1)\widetilde{F}(x_1)t^2,b).
\end{align*}
Hence $t=1$ and from the above
\[\bigl(\forall x\in \widehat{L}^r\setminus\aclk(a)\bigr)\quad F(\aclk(x))=\aclkp(\widetilde{F}(x)).\]
What is left is to show our claim for points from $\aclk(a)\setminus\widehat{K}$. Let $a'\in\aclk(a)\setminus\widehat{K}$. Take independent points $t,s\in\widehat{L}^r\setminus\aclk(a)$, then \[\aclk(a')=\aclk(t,ta')\cap\aclk(s,sa'),\] so as $ta',sa'\in  \widehat{L}^r\setminus\aclk(a)$ from the preceding result we have
\[F(\aclk(a'))=\aclkp(\widetilde{F}(t)),\widetilde{F}(ta'))\cap\aclkp(\widetilde{F}(s), \widetilde{F}(sa'))=\aclkp(\widetilde{F}(a')).\]

The observation that $\GG(L/K)=\GG(\widehat{L}^r/K)$ finishes the proof of (ii).

(iii) It follows immediately from (ii) that $H$ is an epimorphism. Let $f\in\ker H$. Then $j(x,a)=f(j(x,a))=j(f(x),f(a))$, so from \cite[Lemma 2.5]{EH2} there is $n\in\mathbb{Z}$ such that $f(x)=x^{p^n}$ and $f(a)=a^{p^n}$. But $x$ and $a$ were arbitrary (independent), so $f=\text{Frob}^n$.
\end{proof}

\begin{acknowledgements}
I would like to thank Ludomir Newelski for introducing me to this subject and for guidance in preparing this paper.
\end{acknowledgements}

\end{document}